\documentclass[11pt]{article}
\usepackage{amsfonts,amsmath,amssymb,amsthm, amscd}
\usepackage[dvips]{epsfig,graphics}

\numberwithin{equation}{section}

\newtheorem{theorem}{Theorem}[section]
\newtheorem{prop}[theorem]{Proposition}
\newtheorem{lemma}[theorem]{Lemma}
\newtheorem{cor}[theorem]{Corollary}

\theoremstyle{definition}
\newtheorem{definition}[theorem]{Definition}

\newtheorem{remark}[theorem]{Remark}

\def\<{{\langle}}
\def\>{{\rangle}}

\def\e{{\epsilon}}

\def\S{\mathbb S}

\def\L{{\Lambda}}

\def\e{\epsilon}

\def\bs{\bigskip}
\def\Z{\mathbb Z}

\def\D{{\Delta}}

\begin{document}

\title{On modules over Laurent polynomial rings}

\author{Daniel S. Silver \and Susan G. Williams\thanks{Both authors partially supported by NSF grant
DMS-0706798.} \\ {\em
{\small Department of Mathematics and Statistics, University of South Alabama}}}

\maketitle 

\begin{abstract} A finitely generated $\Z[t, t^{-1}]$-module without $\Z$-torsion and having nonzero order $\D(M)$ of degree $d$ is determined by a pair of sub-lattices of $\Z^d$. Their indices are the absolute values of the leading and trailing coefficients of $\D(M)$. This description has applications in knot theory.

  \end{abstract} 

\section{Introduction.} 

The complement of a knot in the 3-sphere fibers over the circle if and only if its universal abelian cover is a product. In that case, the homology of the cover is a finitely generated as an abelian group, and the order of the homology as a $\Z[t, t^{-1}]$-module---the Alexander polynomial of the knot---is monic.

We are motivated by a pair of simple questions: (1) What is the significance of the leading coefficient of the Alexander polynomial of a knot? (2) If $k_1$ and $k_2$ are two knots such that the leading coefficient of the Alexander polynomial of $k_1$ is larger in absolute value than that of $k_2$, in what sense is $k_1$ further away from being fibered than $k_2$? 

An elementary proposition about the structure of finitely generated modules over $\Z[t^{\pm 1}]$ provides answers. The proposition has proved useful in recent study of twisted Alexander invariants \cite{swTAIPR}, but it has not previously appeared in the literature. 

We are grateful to the referee for helpful remarks. 

\section{Main result}  

We denote the ring $\Z[t, t^{-1}]$ of Laurent polynomials by $\L$. It is a Noetherian unique factorization domain. 

Let $M$ be a finitely generated $\L$-module. It can be described as cokernel of a homomorphism 
$$\L^r\  {\buildrel A \over \longrightarrow}\ \L^s,$$
where $A$ is an $r \times s$ {\it presentation matrix} with entries in $\L$. Without loss of generality, we can assume that $r \ge s$. 

\begin{definition} The {\it order} of $M$ is the greatest common divisor of the $s \times s$ minors of the matrix $A$. \end{definition}

We denote the order of $M$ by $\D(M)$. It is well defined up to multiplication by units $\pm t^i$ in $\L$.  The order is nonzero if and only if 
$M$ is a $\L$-torsion module.

In the case that $M$ has a square presentation matrix, the order $\D(M)$ is simply the determinant of the matrix. In general, a finitely generated $\L$-torsion module $M$ admits a square presentation matrix if and only if it has no nonzero finite submodule (see page 132 of \cite{hillman}). 

If $M$ does not have a square presentation matrix, then we can pass to a quotient module that does without affecting the order. To do this, let $zM$ be the maximal finite $\L$-submodule of $M$. Define $\bar M$ to be $M/zM$. One checks easily that $\bar M$ has no nonzero finite submodule. Also, $\D(zM)=1$ (see Theorem 3.5 of \cite{hillman}). The following lemma is well known (see \cite{milnor}). It implies that  
$\D(\bar M) = \D(M)$. 

\begin{lemma} \label{quotient} If $N$ is a submodule of $M$, then 
$$\D(M)=\D(N)\D(M/N).$$ \end{lemma}

Free products with amalgamation are usually defined in the category of groups (see for example \cite{ls}). However, they can also be defined in the simpler category of abelian groups. If $B, B'$ and $U$ are  abelian groups and $f: U \to B, g: U \to B'$ are homomorphisms, then 
we consider the external direct sum $B \oplus B'$ modulo the elements $(f(u), -g(u))$, where $u$ ranges over a set of generators for $U$. We denote the resulting group by $B \oplus_U B'$ when the amalgamating maps $f, g$ are understood.

\begin{lemma}\label{inject} When $g$ is injective, the natural map
$b \mapsto (b, 0)$  embeds $B$ in $B \oplus_U B'$.  Similarly, when $f$ is injective, $b' \mapsto (0, b')$ is an embedding. \end{lemma}

\begin{proof} Assume that $(b, 0)$ is trivial in $B \oplus_U B'$ for some $b \in B$. Then $(b, 0)$ is a linear combination of relators 
$$ (b,0) = \sum_{i=1}^k (f(u_i), -g(u_i)),$$
where $u_1, \ldots, u_k \in U$. The right-hand side can be rewritten
as $(f(u), -g(u))$, where $u = \sum_{i=1}^k u_i$. Hence $g(u)=0$. 
Since $g$ is injective, $u=0$, and so $f(u)=0$ also. 

The proof of the second statement is similar. \end{proof} 

Free products with amalgamation can be formed using infinitely many factors in an analogous fashion. Here we consider groups of the form
$$\cdots \oplus_U B \oplus_U B \oplus_U \cdots$$
with identical amalgamations $B\ {\buildrel
g \over \leftarrow}\  U\ {\buildrel f \over \rightarrow}\ B$.  There is a natural 
automorphism $\mu$ that shifts coordinates one place to the right. It induces the structure of a $\L$-module.

\begin{prop} \label{main} Assume that $M$ is a finitely generated $\L$-module. \begin{itemize}
\item{} There exists a pair $U, B$ of finitely generated abelian groups and monomorphisms $f, g: U \to B$ such that $M$ is isomorphic to the infinite free product 
\begin{equation}\label{decomp}\cdots \oplus_U B \oplus_U B \oplus_U\cdots \end{equation}
with identical amalgamations $B\ {\buildrel
g \over \leftarrow}\  U\ {\buildrel f \over \rightarrow}\ B$;
\item{} If $U$ is generated by $q$ elements, then $\deg \D(M) \le q.$
\end{itemize}  \end{prop}

\begin{cor} \label{cor} Let $M$ be a finitely generated $\L$-module. 
Assume that $\D(M)= c_0 +\cdots + c_d t^d$, with $c_0c_d \ne0$, and ${\rm gcd}(c_0, \ldots, c_d)=1.$  Then  there exist monomorphisms $f, g: \Z^d \to \Z^d$ such that: \begin{itemize}
\item{} $\bar M$ is isomorphic to the infinite free product 
$$\cdots \oplus_{\Z^d} \Z^d \oplus_{\Z^d} \Z^d \oplus_{\Z^d}\cdots$$
with identical injective amalgamations $\Z^d\ {\buildrel
g \over \leftarrow}\ \Z^d\ {\buildrel f \over \rightarrow}\ \Z^d$;
\item{} $\D(M) = \det(t g- f)$;
\item{}
$|c_0|= |\Z^d: f(\Z^d)|$ and  $|c_d|= |\Z^d: g(\Z^d)|$. 
\end{itemize}
\end{cor} 

\begin{remark} It is well known that the coefficients of the Alexander polynomial of a knot are relatively prime and also palindromic. Moreover, the homology of the universal abelian cover of the knot complement is finitely generated as an abelian group if and only $|c_d|\ (= |c_0|) =1$. (See \cite{bz}, for example.)  In such case, the amalgamating maps $f, g$ in the decomposition above identify adjacent factors of $B$. When the Alexander polynomial is not monic, the index $|c_d|$ measures the ``gaps" between the  factors. Corollary \ref{cor} makes the last statement precise. 
\end{remark}



\section{Proof of Proposition \ref{main} and Corollary \ref{cor}.}  The $\L$-module $M$ can be described by generators $a_1, \ldots, a_n$ and relators $r_1, \ldots, r_m$, for some positive integers $n \le m$. For $1\le i\le n$ and 
$\nu \in \Z$, let $a_{i, \nu}$ denote the element $t^\nu a_i$. Similarly, for $1 \le j \le m$, let $r_{j, \nu}$ denote the relator $t^\nu r_j$ expressed in terms of the generators $a_{i, \nu}$. By multiplying at the outset each  $a_i$ by a suitable power of $t$, we can assume without loss of generality that $a_{i, 0}$ occurs in at least one of $r_{1, 0}, \ldots, r_{m, 0}$ but $a_{i, \nu}$ with $\nu < 0$  does not.  Let $\nu_i$ be the largest second index of $a_{i, \nu}$ occurring in any of $r_{1, 0}, \ldots, r_{n, 0}$. 

Define $B$ to be the abelian group with generators $$a_{1, 0}, \dots, a_{1, \nu_1}, \ldots, a_{n, 0}, \ldots, a_{n, \nu_n}$$ and relators $r_{1, 0}, \ldots, r_{m, 0}$.  Let $U$ be the free abelian group generated by $$a_{1, 0}, \dots, a_{1, \nu_1-1}, \ldots, a_{n, 0}, \ldots a_{n, \nu_n-1}.$$
(If some $\nu_i =0$, then all $a_{i, \nu}$ are omitted.)
Define a homomorphism $f: U \to B$ by mapping each $a_{i, \nu} \mapsto a_{i, \nu}$. Define $g: U \to B$ by $a_{i,\nu} \mapsto a_{i, \nu+1}$.

If either $f$ or $g$ is not injective, then we apply the following {\it reduction operation}: replace $U$ by the quotient abelian group  $U/({\rm ker} f + {\rm ker}\ g)$. Replace $B$ by $B/(f({\rm ker}\ g)+ g({\ker}\ f))$, and $f, g$ by the unique induced homomorphisms. If again $f$ or $g$ fails to be injective, repeat this operation. By the Noetherian condition for abelian groups, we obtain injective maps $f, g$ after finitely many iterations.

Consider the free product with amalgamation
\begin{equation}\label{decomposition}\cdots \oplus_U B \oplus_U B \oplus_U \cdots.\end{equation}
The amalgamation maps $B\ {\buildrel
g \over \longleftarrow}\ U\ {\buildrel f \over \longrightarrow}\ B $ are injective as a consequence of the reduction operation.  

A surjection $h$ from $M$ to the module described by (\ref{decomposition}) is easily defined:
Map each generator $a_{1, 0}, \dots, a_{1, \nu_1}, \ldots, a_{n, 0}, \ldots, a_{n, \nu_n}$ to the generator with the same name in a chosen factor $B$. For any $\nu \in \Z$, map $t^\nu a_{1, 0}, \dots,\break  t^\nu a_{1, \nu_1}, \ldots, t^\nu a_{n, 0}, \ldots, t^\nu a_{n, \nu_n}$ to the same generators but in the shifted factor $\mu^\nu(B)$. The amalgamating relations of (\ref{decomposition}) ensure that $h$ is well defined. The relations that result from the reduction operation are mapped trivially, since the image under $f$ (resp.\  $g$) of any element in the kernel of $g$ (resp.\  $f$) is trivial in $M$. Hence $h$ is an isomorphism. The first part of Proposition \ref{main} is proved.

Assume that $U$ is generated by $q$ elements while $B$ is described by $Q$ generators and $P$ relations. We take $q \le Q$. A $(P+q) \times Q$ presentation matrix for the $\L$-module $M$ is easily constructed. 
Columns correspond to the generators of $B$. Rows correspond to the relations of $B$ together with relations of the form $t g(u) = f(u)$, where $u$ ranges over generators of $U$.   Without loss of generality, we can assume that $P + q \ge Q$. The order $\D(M)$ is the greatest common divisor of the $Q \times Q$ minors. The greatest power of $t$ that can occur is $q$. No negative powers arise. Hence $\deg \D(M) \le q$.

In order to prove Corollary \ref{cor}, we use Proposition \ref{main} to write
$\bar M = \cdots \oplus_U B \oplus_U B \oplus_U \cdots$ with identical injective amalgamating maps. 
By Lemma \ref{inject}, $B$ is a subgroup of $\bar M$. Since the coefficients of $\D(\bar M)$ have no nontrivial common factor, $\bar M$ is $\Z$-torsion free by \cite{crowell}. Hence  $B$ and $U$ are  finitely generated free abelian groups. Moreover, the rank of $U$ must equal that of $B$, since otherwise $\D(M)$ would vanish. The amalgamations $B\ {\buildrel
g \over \leftarrow}\ U\ {\buildrel f \over \rightarrow}\ B $ induce a defining set of module relators $t g(u) - f(u)$, where $u$ ranges over a basis for $U$. Hence the order of $\bar M$, which is equal to $\D(M)$, is $\det(t g - f)$. The remaining statement of Corollary \ref{cor} follows immediately.

\section{Examples.} 

\indent\indent


Some hypothesis on the finitely generated $\L$-module $M$ is needed for the conclusion of Corollary \ref{cor} to hold, as the next example demonstrates. \bs

(1) Consider the cyclic $\L$-module $M= \< a \mid 2a =0\>.$ Since $M$ has a square presentation matrix, $zM$ is trivial and hence $\bar M = M$. In this case, $B \cong \Z/2\Z$ and $U$ is trivial. Since every nonzero element of $M$ has order $2$, it contains no free abelian subgroup $B$. Lemma \ref{inject} shows that the conclusion of Corollary \ref{cor} fails in this example. \\

The hypothesis of Corollary \ref{cor} on the coefficients of $\D(M)$ is not a necessary condition, as we see in the next example. \\

(2) Consider the cyclic $\L$-module $M= \< a \mid (2t-2)a=0 \>.$
As in the previous example, $zM = 0$ and hence $\bar M = M$. 
Although the coefficients of $\D(M)= 2t-2$ are not relatively prime, 
$M$ can be expressed in the form (\ref{decomp}) with $B$ and $U$ infinite cyclic  groups and $f, g: U \to B$ given by $u \mapsto 2 u$. \\

(3) Let $G$ be a finitely presented group and $\e: G \to \Z$ an epimorphism. Denote the kernel of $\e$ by $K$. Choose an element $x$ such that  $\e(x)=1$. The abelianization $M=K/K'$ is a finitely generated $\L$-module with
$t (a+K') = xax^{-1}+K'$ for all cosets $a+K'$. The module depends only on $G, \e$ and not on the choice of $x$. 

Consider the case that $G=\pi_1(\S^3 \setminus k)$ is the group of an oriented knot $k$, the map $\e$ is abelianization, and $x$ is a meridian. 
The module $M$ is the homology of the infinite cyclic cover of $
\S^3 \setminus k$. It is well known that $M$ is $\Z$-torsion free \cite{crowell2}.  In particular, the conclusion of Corollary \ref{cor} holds. 

Splitting the 3-sphere along a Seifert surface $S$ for the knot with minimal genus $g$ produces a relative cobordism $Y$ between two copies $S_-$, $S_+$ of $S$.  Let $B = H_1 Y$ and $U = H_1 S$. By Alexander duality, $U \cong B \cong \Z^{2g}$. Let $f: H_1 S_- \to B$ and $g: H_1 S_+ \to B$ be homomorphisms induced by inclusion. Then $M$ can be expressed in the form (\ref{decomp}). However, $f, g$ are generally not injective. The reduction operation described in the proof of Proposition \ref{main} can reduce the rank of $B$. In fact, when $k$ has Alexander polynomial $1$, the group $B$ must become trivial.


\end{document}